\documentclass[leqno,11pt]{amsart}
\usepackage{latexsym,amsfonts,amsmath}

\newtheorem{theorem}{Theorem}
\newtheorem{cor}[theorem]{Corollary}
\newtheorem{conj}[theorem]{Conjecture}
\newtheorem{lemma}[theorem]{Lemma}
\newtheorem{prop}[theorem]{Proposition}
\theoremstyle{definition}
\newtheorem{definition}{Definition}[section]
\theoremstyle{remark}
\newtheorem{remark}{Remark}[section]
\numberwithin{equation}{section}
\newcommand{\nc}{\newcommand}
 \newcommand{\R}{{\mathbb R}}

\nc{\ar}{\rightarrow}
\nc{\hd}[1]{\mathcal{H}(#1)}
\nc{\kd}[1]{\mathcal{K}(#1)}
\nc{\ds}{d_{S'}}
\nc{\sumj}{\sum_{j=1}^n}
\nc{\B}{\mathcal{B}}
\nc{\W}{\mathcal{W}}
\nc{\fb}{S\mathcal{B}}
\nc{\fw}{S\mathcal{W}}
\nc{\tc}{\widetilde{C}}
\nc{\ws}{\widetilde{S}}
\nc{\ot}{\tilde\omega}
\title[Exceptional points for Lebesgue's density theorem]{Exceptional points for Lebesgue's density theorem on the real line}
  \author{Andr{\'a}s Szenes} 
\address{BME Mathematics Institute, Budapest, Hungary; email:
szenes@math.bme.hu}
\thanks{ The support of the Hungarian National Science
Foundation Grant OTKA-46365 is acknowledged.}
\begin{document}
\maketitle

\setcounter{section}{-1}
\section{Introduction and notation}
\label{intro}

\subsection{Formulation of the problem}
Denote by $\lambda$ the Lebesgue measure on the real line. We will
call a measurable set $S\subset\R$ {\em nontrivial} if neither $S$ nor
$\R\setminus S$ is of measure zero. A point $p\in\R$ is called a {\em
  density point} of $S$ if
\[
\lim_{\epsilon\ar0}\frac{\lambda(I_\epsilon(p)\cap S)}{2\epsilon}=1,
\]
where $I_\epsilon(p)$ is the interval $(p-\epsilon,p+\epsilon)$.

The well-known Lebesgue density theorem, in a somewhat weakened form,
states that
\smallskip

{\em   For any measurable set $S\subset\R$, almost all points $p\in\R$ are
  either density points of $S$ or density points of $\R\setminus S$.}

\smallskip It is a natural problem to investigate the set of what we
will call {\em exceptional points} for $S$, i.e. points which are
neither density points of $S$, nor those of $\R\setminus S$. Note that
this is a topological notion, since as far as measure theory is
concerned, there are no such exceptional points.

First, we quantify the notion of exceptional point:
given a measurable $S\subset\R$ and $0\leq \delta\leq1/2$, we will
call $p\in \R$ a $\delta${\em-exceptional point} for $S$ if
\[
\delta\leq\liminf_{\epsilon\ar0}\frac{\lambda(I_\epsilon(p)\cap S)}
{2\epsilon}\leq\limsup_{\epsilon\ar0}\frac{\lambda(I_\epsilon(p)
\cap S)}{2\epsilon}\leq1-\delta. 
\]

Let $0\leq\delta\leq1/2$. In this article, we will be studying the
statement
\[ \hd\delta:\;\text{\bf There is a }\delta\text{\bf -exceptional
point for every nontrivial }S\subset\R.
\]
Clearly, if $\delta_1>\delta_2$ then $\hd{\delta_1}$ implies
$\hd{\delta_2}$. The central problem we are addressing is finding the
universal constant $\delta_{\mathcal{H}}$:
\[ \delta_{\mathcal{H}}=\sup\{\delta|\; \hd\delta\text{ is true}\}.\]

\subsection{The history of the problem} The problem of determining the
constant $\delta_{\mathcal{H}}$ was introduced and studied in
\cite[\S4]{vik}; in this paper Victor Kolyada showed that
\[
1/4\leq\delta_{\mathcal{H}}\leq(\sqrt{17}-3)/4\sim0.2807764
\]
On his suggestion, the question of proving the inequality
$1/4\leq\delta_{\mathcal{H}}$ became one of the problems in the 1983
Schweitzer competition (cf. \cite[Problem 9, 1983]{schw}), a contest
for mathematics undergraduates in Hungary. As it turned out, the
author could not solve this problem at the time, and, as a result,
failed to win the first prize in the competition. Probably, to some
extent motivated by this disappointment, the author undertook a
thorough study of the problem after the competition, and this led to
the result obtained in 1984, which, with apologies for the
considerable delay, we submit in the present paper.

\subsection{Results, and contents of the paper} There is a simple
analytic proof of the fact that $\delta_{\mathcal{H}}\geq 1/4$; we
recall this proof in \S1. In \S2 we describe a combinatorial
restatement of our problem, and using this combinatorial approach, in
\S3, we give an upper bound on $\delta_{\mathcal{H}}$. We conjecture
that this upper bound, which is a solution of a cubic equation, and is
approximately 0.272, is, in fact, the value of
$\delta_{\mathcal{H}}$. The main result of the paper is described in
the last section, where we prove a lower bound on
$\delta_{\mathcal{H}}$. This lower bound is also a solution of a cubic
equation; its value is about 0.263.

{\sc Notation and conventions}: In this article, every set is assumed
to be measurable. All intervals will be considered open. The length of
an interval $J$ will be denoted by $|J|$. We denote by $I_\epsilon(p)$
the $\epsilon$-neighborhood of the point $p\in\R$, i.e.  the interval
$(p-\epsilon, p+\epsilon)$.  Given an interval $I\subset\R$ and a
subset $H\subset\R$, denote by $\lambda(H|I)$ the relative measure of
$H$ in I, i.e.
\[ \lambda(H|I) = \frac{\lambda(H\cap I)}{|I|}.
\]
Given a set $S\subset\R$ and a number $a\in\R$ we denote by $a+S$ the
set $\{a+x;\;x\in S\}$ and by $a-S$ the set $\{a-x;\;x\in S\}$.

{\sc Acknowledgment.} We would like to thank Victor Kolyada for useful
comments and references, and extend our gratitude to Mikl\'os
Laczkovich for his help and encouragement.  

\section{The solution of the Schweitzer problem}

\begin{prop} \label{quarter}
  The statement $\hd{1/4}$ is true.
  \end{prop}
  Let us see the proof. We are given a nontrivial $S\subset\R$, and we are
  looking for a 1/4-exceptional point for $S$.  Let $a$ be a density
  point for $S$ and $b$ be a density point for the complement of $S$.
  Without loss of generality we may assume that $a=0$ and $b=1$.
  Denote by $\ws$ the truncated set $\ws=(-\infty,0)\cup
  S\setminus(1,\infty)$ and let $d_{\ws}(x) =
  \lambda(\ws\cap(x,\infty))$.  The function
\[f(x)=d_{\ws}(x)+x/2\]
goes to infinity linearly as $x\ar\pm\infty$, and its derivative is
negative at 0, and positive at 1. This implies that $f(x)$ has
a global minimum at a point $p$ in the {\em interior} of the interval
$(0,1)$. Now, given an arbitrary $\epsilon>0$, we have
\begin{multline}
\lambda(\ws|I_\epsilon(p))\geq\frac{\lambda((p-\epsilon,p)
\cap \ws)}{2\epsilon}=\frac{d_{\ws}(p-\epsilon)-d_{\ws}(p)}{2\epsilon}
\\ =\frac{(d_{\ws}(p-\epsilon)+(p-\epsilon)/2)-(d_{\ws}(p)+p/2)}{2\epsilon}
+\frac14=\frac{f(p-\epsilon)-f(p)}{2\epsilon}+\frac14\geq\frac14;
\end{multline}
similarly, one sees that
\[ \lambda(\ws|I_\epsilon(p)) \leq\frac34.\]
As $0<p<1$, the sets $S$ and $\ws$ coincide near $p$, and thus $p$ is
a $1/4$-exceptional point for $S$. This proves that $\hd{1/4}$
holds.\qed

It does not appear that this proof can be improved upon easily, thus
it seems natural to conjecture that, in fact,
$\delta_{\mathcal{H}}=1/4$.  Thus we were very surprised to discover
otherwise. To explain the reasons behind this phenomenon, we first
recast the problem in a discrete form.

\section{Combinatorial restatement}
Based on an idea of Mikl\'os Laczkovich, we formulate a combinatorial
problem, which turns out to be equivalent to determining whether
$\hd\delta$ is true (also cf. \cite[\S4]{vik}).

Given a finite, increasing sequence of positive real numbers, 
\[ 0<a_1<b_1<\dots<a_r<b_r,\]
we call the union of intervals 
\[ C=(-\infty,0)\cup\bigcup_{i=1}^r(a_i,b_i) \]
a {\em configuration}, and the elements of the sequence, including 0,
the {\em endpoints} of $C$.

Given $\delta$, $0\leq\delta\leq1/2$, we denote by $\kd \delta$ the
following statement:
\\
$\kd\delta$:\;
{\bf For every configuration $C$, there is an endpoint $c$ such that
\[ \delta\leq \lambda(C|I_\omega(c))\leq1-\delta\text{ for all } \omega>0.\]
}

\smallskip 
For the convenience of the reader, we write down the {\bf opposite} of
$\kd\delta$ as well:\\
{\em  There exists a configuration $C$ such that for
  every endpoint $c$ of $C$ there is a positive radius $\omega(c)$
  such that $\lambda(C|I_{\omega(c)}(c))\notin[\delta,1-\delta]$.}

\smallskip
Again, clearly $\kd{\delta_1}$ implies $\kd{\delta_2}$ if
$\delta_1>\delta_2$. Set
$\delta_{\mathcal{K}}=\sup\{\delta>0;\,\kd{\delta}\text{ true}\}$.
\begin{prop}\label{hequalk}
  We have $\delta_{\mathcal{H}}=\delta_{\mathcal{K}}$.
\end{prop}
\begin{proof}
  First we show that if $\hd \delta$ is false, then so is
  $\kd{\delta+\tau}$ for any $\tau>0$. Assume that $S$ is a
  counterexample to $\hd{\delta}$. Using the cut-off construction at
  the beginning of Proposition \ref{quarter}, without loss of
  generality, we can assume that $(1,\infty)\cap S=\emptyset$ and
  $(-\infty,0)\subset S$.  Then for every $x$ in the closed interval
  $[0,1]$, there exists a radius $\epsilon(x)$ such that
  $\lambda(S|I_{\epsilon(x)}(x))\notin[\delta,1-\delta]$.  At the cost
  of increasing $\delta$, one may put a uniform lower bound on
  $\epsilon(x)$. Indeed, fix a small $t>0$. It is easy to check that
  for $y\in I_{t\epsilon(x)}(x)$ we have
  $\lambda(S|I_{\epsilon(x)}(y))\notin[\delta+t,1-\delta-t]$. Since
  $[0,1]$ is compact, it is covered by finitely many of the intervals
  $I_{t\epsilon(x)}(x)$. Pick such a finite cover and denote by $\eta$
  the least of the radii $\epsilon(x)$ in it. Then for each $y\in
  [0,1]$ there is an $x\in[0,1]$ such that $y\in t\epsilon(x)$,
  $\epsilon(x)\geq\eta$  and 
\[\lambda(S|I_{\epsilon(x)}(y))\notin[\delta+t,1-\delta-t].
\] Finally, by approximating $S$ with a finite union of intervals, we
  can find a configuration $C$ such that for any interval $I$ we have
  \[ |\lambda(I\cap C)-\lambda(I\cap S)|<t\eta.  \] Then by applying
  to each endpoint of $C$ the last two inequalities, we can convince
  ourselves that  $C$ provides a counterexample to $\kd{\delta+2t}$.
This clearly shows that $\delta_\mathcal{K} \leq \delta_\mathcal{H}$.

  Now we prove the opposite inequality. Assume that the configuration $C$ is a
  counterexample to $\kd \delta$. This means that for each endpoint
  $c$ of $C$ there is a radius $\omega(c)>0$ such that
  $\lambda(C|I_{\omega(c)}(c))\notin[\delta,1-\delta]$). Denote the
  least and greatest among the positive numbers $\omega(c)$ by
  $\omega_{\min}$ and $\omega_{\max}$ respectively.

  Without loss of generality, we can assume that
  $C\subset(-\infty,1)$; let $\tilde C=C\cap(0,1)$. Fix a small
  $\epsilon>0$ and let $H_1=\tilde C$. We define a finite disjoint
  union of intervals $H_n$ by induction as follows: write
  $H_n=\cup_{j=1}^{r(n)}(a_j(n),b_j(n))$ and let
\[  H_{n+1} = \bigcup_{j=1}^{r(n)}\left([a_j(n)-\epsilon^n\tilde C]
\cup(a_j(n),b_j(n))\cup[b_j(n)+\epsilon^n\tilde C]\right).
\]
In particular, $H_n\subset H_{n+1}$. 

Finally, let $H=\cup_{n=1}^\infty H_n$. We will now show that for any
$\tau>0$ one can choose a sufficiently small $\epsilon>0$ such that
$H=H(\epsilon)$ is a counterexample to $\hd {\delta+\tau}$. Pick an
arbitrary point $x\in\R$. We need to compute $\liminf/\limsup$ of the
density of the set $H$ around $x$.  Clearly, we can assume that $x$ is
a boundary point of $H$, otherwise the density is 0 or 1.

Pick a positive integer $n$ and denote by $v=v_n$ the endpoint of
$H_n$ closest to $x$. Since $C$ is a counterexample to $\kd \delta$,
there is a a radius $\omega=\omega_n$,
$\omega_{\min}\leq\omega\leq\omega_{\max}$ such that
\begin{equation}
  \label{hn}
\lambda(H_n|I_{\epsilon^{n-1}\omega}(v))<\delta.  
\end{equation}
  For simplicity of notation, we will suppress the other possibility:
  $>1-\delta$. We would like to estimate $\lambda(H|I_{\epsilon^{n-1}
    \omega}(x))$. 
  
  First, using the trivial bound $\lambda(\tilde C)\leq1$, we obtain
  \begin{equation}
    \label{hhn}
\lambda((H\setminus H_n)\cap I_{\epsilon^{n-1}\omega}(v)) 
 <\epsilon^{n-1} \frac{ M\epsilon}{1-M\epsilon},    
  \end{equation}
where $M$ is the number of endpoints of $C$.

Next, we can estimate the distance between $x$ and $v$ as
\begin{equation}
  \label{xv}  |x-v|\leq \epsilon^n.
\end{equation}
Combining the inequalities \eqref{hn}, \eqref{hhn} and \eqref{xv}, a
short computation shows that 
\[
\lambda(H|I_{\epsilon^{n-1}\omega}(x))<\delta+
\frac\epsilon{2\omega}\left(1+\frac M{1-M\epsilon}\right).
\]
Thus given any $\tau>0$,  we can choose $\epsilon$ sufficiently small,
so that we have
\[\lambda(H|I_{\epsilon^{n-1}\omega_n}(x))\notin(\delta+\tau,1-\delta-\tau)\]
for the sequence of intervals constructed above. Since clearly
$\epsilon^{n-1}\omega_n\rightarrow0$, we can conclude that
$\delta_\mathcal{K} \geq \delta_\mathcal{H}$, and this completes the
proof.
\end{proof}

\section{An upper bound}\label{upper}

The main goal of this article is to estimate the constant
$\delta_\mathcal{H}$ introduced in \S\ref{intro}. The rather
``natural'' proof of Proposition \ref{quarter} seems to suggest that 
$\delta_\mathcal{H}=1/4$. In the next section, we will prove,
however, that $\delta_\mathcal{H}>1/4$!

Proposition \ref{hequalk} shows that we can study the constant
$\delta_\mathcal{K}$ instead of $\delta_\mathcal{H}$.  The following
statement provides an upper bound for $\delta_\mathcal{K}$.
\begin{prop}
  \label{counter}
If $(2\delta)^3+(2\delta)^2+2\delta>1$, then there is a counterexample
to $\kd \delta$.
\end{prop}
\begin{remark}
  This provides the bound $\delta_\mathcal{K}<0.2719$.
\end{remark}
\begin{proof}
  We construct a configuration $C(m,s,N)\subset(-\infty,1)$ depending
  on 2 parameters, $0<m,s<1$ and a large integer $N$. The construction
  goes as follows. We consider the interval $(1-m,1)$, and divide it into
  $N$ equal parts. Next we break each of these parts into two: an initial piece
  proportional to $s$ and a final piece, proportional to $1-s$, and
  then take the union of these initial pieces:
\[ C(m,s,N)\setminus(-\infty,0)=
\left\{x\in(1-m,1);\;0<\left\{\frac{N(x+m-1)}{m}\right\}<s\right\},
\]
where $\{y\}$ stands for the fractional part of the  real number $y$.
Then we can compile the following table: the first column lists the
endpoints of $C(m,s,N)$, the second  a certain chosen radius, and the last
one twice the corresponding density.
\begin{center}
\begin{tabular}{l|l|l}
endpoint $v$ & radius $r$ & $2\lambda(C(m,s,N)|I_r(v))$\\ \hline
0 & 1 & $sm+1$\\
$1-m$ & $m$ & $2-(1/m-s)$\\
$\sim1$ & 1 & $\sim sm$\\
all other & $sm/N$ & $2-(1/s-1)$
\end{tabular}
\end{center}
The third line of the table represents the last endpoint of
$C(m,s,N)$; it approaches 1 as $N\rightarrow\infty$ and the
corresponding density has been computed in this limit as well.
It is clear that all but this endpoint give densities $>1/2$, and that
the first density: $sm+1$, is always greater than the second:
$2-(1/m-s)$. 

Then a simple argument shows that the optimal configuration (in the
limit when $N\rightarrow\infty$) is achieved when 
\begin{equation}
  \label{triple}
\frac1m-s=sm=\frac1s-1.  
\end{equation}
Indeed, it is sufficient to check that the gradients of the three
two-variable functions which appear here are never
collinear. Eliminating $m$ from \eqref{triple} we obtain
\[ 2s^3-2s^2+2s=1.\] This quickly leads to the equation
\[ q^3+q^2+q=1\]
for the parameter $q=1/s-1$, which represents twice the density. This
completes the proof.
\end{proof}
We conjecture that this is, in fact, an optimal construction.
\begin{conj}
  \label{conj}
  The universal constant $\delta_\mathcal{K}$ is the only real root of
  the cubic equation
\[ (2\delta)^3+(2\delta)^2+2\delta=1.\]
\end{conj}
We have not been able to prove this conjecture; see, however, Remark
\ref{last}. 

\section{The Main Result}

\begin{theorem}
$\kd \delta$ is true if $4\delta^3+2\delta^2+3\delta<1$.
\end{theorem}
\begin{remark} The theorem provides the lower bound
  $\delta_\mathcal{K}>0.2629$. 
\end{remark}
We start with a simple Lemma.
\begin{lemma}
  \label{first}
Suppose that an interval $I$ is represented as a not necessarily
disjoint union of intervals: $I=\cup_{j=1}^nI_j$. Assume that $0<\delta<1$,
and let $B$ be a measurable set such that $\lambda(B|I_j)\geq
1-\delta$ for $j=1,\dots,n$. Then 
\[ \lambda(B|I)\geq \frac{1-\delta}{1+\delta}.
\]
\end{lemma}
\begin{proof}
  Without loss of generality we can assume that $I=(0,1)$, and that
  our system of intervals $I_j=(a_j,b_j)$, $j=1,\dots,n$, satisfies
\begin{enumerate}
\item $a_j<a_{j+1}$, for $j=1,\dots,n-1$, i.e. the left endpoints form
  an increasing sequence, and
\item $I_j\cap I_{j+2}=\emptyset$ for $j=1,\dots,n-2$.
\end{enumerate}
Indeed, the first condition can be satisfied by renumbering the
intervals, and the second by eliminating intervals which are contained
in the union of the rest of the system. Introduce the following
parameters of the system: setting $I_0=I_{n+1}=\emptyset$, for $1\leq
j\leq n$ let
\begin{eqnarray*}
& x_j=\lambda(I_j\cap I_{j+1}),&x_j^B=\lambda(I_j\cap I_{j+1}\cap
B),\\
 & y_j = \lambda(I_j\setminus(I_{j-1}\cup I_{j+1})),&
y^B_j = \lambda((B\cap I_j)\setminus(I_{j-1}\cup I_{j+1})).
\end{eqnarray*}
Using these parameters, we can rewrite the inequality
$\lambda(B|I_j)\geq1-\delta$ as 
\[ x^B_{j-1}+y^B_j+x^B_j\geq(1-\delta)(x_{j-1}+y_j+x_j).
\]
Summing these inequalities for $j=1,\dots,n$, we obtain
\[ \frac{2x^B+y^B}{2x+y}\geq1-\delta,
\]
where
\[ x=\sumj x_j,\,y=\sumj y_j,\, x^B=\sumj x_j^B,\,y^B=\sumj y_j^B.
\]
Now using the fact that $x+y=1$, and that $x^B\leq x$, we can conclude
that 
\[ \frac{2x^B+y^B}{1+x^B}\geq1-\delta.
\]
Hence
\[ (1+\delta)x^B+y^B\geq 1-\delta,
\]
which implies that 
\[ x^B+y^B\geq\frac{1-\delta}{1+\delta}.
\] This last inequality is exactly the statement of the Lemma.
\end{proof}

Now we begin the proof of the Theorem. Assume that $\kd \delta$ does not
hold for some $0<\delta<\frac12$. Our results so far show that in this
case $1/4< \delta$. Then let
\[  C=(-\infty,0)\cup(a_1,b_1)\cup\dots\cup(a_r,b_r=1)
\]  be a configuration which is a counterexample to $\kd \delta$ with the
{\em least} possible number $r$ of intervals in it. For each endpoint
$p$ of $C$, introduce the set
\[ D_p = \{\omega\in \R_{\geq0};\;
\lambda(C|I_\omega(p))\notin(\delta,1-\delta)\}, 
\]
and let $\omega(p)=\sup D_p$. Note that, by our assumption, $D_p$ is
nonempty for every endpoint $p$ of $C$.
\begin{definition}
  We will call an endpoint $p$ {\em black} if
  $\lambda(C|I_{\omega(p)}(p))\geq1-\delta$, and {\em white} if
  $\lambda(C|I_{\omega(p)}(p))\leq \delta$. Denote the set of black endpoints by
  $\B=\B(C)$, and the set of white endpoints by $\W=\W(C)$.
\end{definition}
Notice that $0$ is a black, while $1$ is a white endpoint.
\begin{lemma}
  \label{second}
If $p$ is a black endpoint and $p\leq1/2$, then either $\omega(p)<p$ or
$\omega(p)\geq1-p$. Similarly, for $p\in\W$ and $p\geq1/2$, we have
$\omega(p)< 1-p$ or $\omega(p)\geq p$.
\end{lemma}
\begin{proof}
Assume that contrary to the statement of the Lemma, there is a
$p\in\B$ such that $\omega(p)\geq p$ and $p+\omega(p)<1$. We will arrive
at a contradiction from these assumptions.

First we observe that we must have $b_i\leq p+\omega(p)\leq a_{i+1}$ for some
$i<r$. Indeed, if $p+\omega(p)$ were an interior point of an interval in $C$,
then for a sufficiently small $\epsilon>0$, the density
$\lambda(C|I_{\omega(p)+\epsilon}(p))$ would be strictly greater than the
density $\lambda(C|I_{\omega(p)}(p))$; this contradicts the definition of
$\omega(p)$ as the maximal radius $\omega$ for which
$\lambda(C|I_{\omega}(p))\geq1-\delta$.

Now we claim that the configuration 
\[ C_{p+\omega(p)} = C\setminus(p+\omega(p),\infty) \]
is a counterexample to $\kd \delta$. For every vertex v of $C_{p+\omega(p)}$, we
need to find an appropriate radius $\ot(v)$, such that
\begin{equation}
  \label{toprove}
\lambda(C_{p+\omega(p)}|I_{\ot(v)}(v))\notin[\delta,1-\delta].  
\end{equation}
It follows from our observation above that the vertices of
$C_{p+\omega(p)}$ form a subset of the vertices of $C$. If $v\in\W(C)$, or
$v\in\B(C)$ and $v+\omega(v)\leq p+\omega(p)$, then then \eqref{toprove} is easy
to satisfy: one chooses $\ot(v)=\omega(v)$. Pick a black vertex $v\in\B(C)$
with $v+\omega(v)\leq p+\omega(p)$. To show that $C_{p+\omega(p)}$ is a
counterexample to $\kd \delta$ we prove that
\[ \lambda(C_{p+\omega(p)}|I_{p+\omega(p)-v}(v))>1-\delta.\]
Indeed, the definition of $\omega(p)$ implies that
$\lambda(C|(p+\omega(p),v+\omega(v)))<1-2\delta$. This, in turn, means that
\[ 1-\delta = \lambda(C|I_{\omega(v)}(v))<\lambda(C|I_{p+\omega(p)-v}(v)).\]

Now observe that the configuration $C_{p+\omega(p)}$ has fewer elements
than $C$. The fact that it provides a counterexample to $\kd \delta$
contradicts  $C$ being  a counterexample with the fewest
possible number of intervals in it. This completes the proof of the Lemma.
\end{proof}

We can divide the set $\{v\in\B; v\leq\frac12\}$ into two
groups: in the first group we collect the endpoints which satisfy
$\omega(v)<v$; the second group will contain the endpoints for which
$\omega(v)\geq v$, in which case $\omega(v)\geq1-v$ according to Lemma
\ref{second}. This second group is always nonempty since 0 is in it.
Introduce a special notation for the largest endpoint from the second
group:
\[ v_\B=\max\{v\in\B;\; v\leq1/2,\,\omega(v)\geq1-v\},\]
and also let
\[ v_\W=\min\{v\in\W;\; v\geq1/2,\, \omega(v)\leq v\}.\]
In addition, set $\rho=\lambda(C\cap(0,1))$ and
$I_\circ=(v_\B,v_\W)$. 
\begin{lemma}
\label{third}
  In the notation introduced above, we have
\[ \frac{1-\rho}{2(1-v_\B)}\leq
\delta\quad\text{and}\quad\frac\rho{2v_\W}\leq \delta.\]
\end{lemma}
\begin{proof}
  It is easy to see that if for a black endpoint $v$ between $0$ and
  $1/2$ we have $\omega(v)\geq1-v$, then
  $\lambda(C|I_{1-v}(v))\geq1-\delta$. This implies the first
  equality. The second one is proved similarly.
\end{proof}
The following statement is the heart of our argument. Its proof will
take up most of the remainder of the paper.
\begin{prop}
  \label{fourth}
\[ \rho\geq\frac{1-\delta}{1+\delta}|I_\circ|\quad\text{or}\quad
\rho\leq\left(1-\frac{1-\delta}{1+\delta}\right)|I_\circ|. \]
\end{prop}
\begin{proof}  
If $C$ has no endpoints inside $I_\circ$, then the statement of
the Proposition is satisfied trivially. We can thus assume that
the set $F$ of endpoints of $C$ inside $I_\circ$ is non-empty:
\[ F = \{v\in\B\cup\W;\;v_\B<v<v_\W\}\neq\emptyset.\]

Now for $v\in\B$ denote by $\mu(v)$ the radius of the interval around
$v$ in which the density of $C$ is maximal. Thus for any $\omega>0$, we
have 
\[ \lambda(C|I_{\mu(v)(v)})\geq \lambda(C|I_{\omega}(v)).\]
Similarly, for $v\in\W$, we denote by $\mu(v)$ the radius of the
interval around $v$ in which the density of $C$ is minimal. 

\begin{lemma}\label{cimp}
  If $p\in F$, then $I_{\mu(p)}(p)\subset(0,1)$.
\end{lemma}
\begin{proof}
Assume that $p\leq\frac12$. Then if $p\in\B$, then $\mu(p)\leq p$
because of the definition of $v_\B$. If $p\in\W$ and
$\lambda(C|I_p(p))\leq\frac12$, then 
$\lambda(C|I_{\omega}(p))$ will increase with $\omega$ for $\omega>p$. This
implies that in this case, again, $\mu(p)\leq p$. The proof in the
case when $p>\frac12$ is analogous.
\end{proof}

Now we construct two subsets $\fb$ and $\fw$ of the interval $(0,1)$
as follows. Let
\begin{eqnarray*}
 &\fb_1 = \cup\{I_{\mu(p)}(p);\;p\in F\cap\B\}, &
\fw_1 = \cup\{I_{\mu(p)}(p);\;p\in F\cap\W\}\\
&\fb_2 = \cup\{(a_i,b_i);\;\fb_1\cap(a_i,b_i)\neq\emptyset\},
&\fw_2 = \cup\{(b_i,a_{i+1});\;\fw_1\cap(b_i,a_{i+1})\neq\emptyset\}\\
&\fb=\fb_1\cup\fb_2,&\fw=\fw_1\cup\fw_2.
\end{eqnarray*}
Clearly, all these sets are unions of intervals.
\begin{lemma}\label{icircin}
  \[ I_\circ\subset \fb\cup\fw\subset(0,1).\]
\end{lemma}
\begin{proof}
  The fact that $\fb,\fw\subset(0,1)$ easily follows from Lemma
  \ref{cimp}. Now let $(a_i,b_i)\subset I_\circ$. Then either $a_i$ or
  $b_i$ is an element of $F$, i.e. lies in the interior of
  $I_\circ$. Assume that $a_i\in F$. If $a_i\in\W$, then 
  $(a_i,b_i)\subset I_{\mu(a_i)}(a_i)$, and thus
  $(a_i,b_i)\subset\fw_1$. On the other hand, if $a_i\in\B$, then
  obviously $(a_i,b_i)\subset \fb_2$. The other case, $b_i\in F$ is
  similar. It is not hard to see that the same method of proof works
  for the intervals of the form $(b_j,a_{j+1})$. This completes the
  proof of the Lemma. 
\end{proof}

\begin{lemma}\label{mainlemma}
  \begin{enumerate}
  \item The set $\fb$ is a union of intervals of the form $(a_i,b_j)$,
  $i\leq j$, while the set $\fw$ is a union of intervals of the form
  $(b_i,a_j)$, $j<i$.
\item 
  Let the intervals $J_\B$ and $J_\W$ be connected components of the
  sets $\fb$ and $\fw$, respectively. Then exactly one of the following 3
  possibilities takes place:
\[ J_\B\cap J_\W=\emptyset\;\text{ or }\;J_\B\subset J_\W\;
\text{ or }\;J_\W\subset J_\B.\] 
 \end{enumerate}
\end{lemma}
\begin{proof}
  To prove the first statement, observe that for $p\in F\cap\B$, the interval
  $I_{\mu(p)}(p)$  has to have its two boundary
  points in the closure of  $C$ in order to conform with the
  definition of $\mu(p)$. These two intervals are subsets of $\fb_2$ by
  construction, and this completes the proof for $\fb$. The proof is
  similar for $\fw$.
  
  Now we turn to the second statement, which is the key to our whole
  argument. It follows from (1) that $J_\B=(a_i,b_j)$ and
  $J_\W=(b_k,a_l)$ for some indices $0\leq i,j,k,l\leq n$.  If the two
  intervals, $J_\B$ and $J_\W$ were not situated as described in the
  statement, then we would have the following two remaining
  possibilities:
  \begin{equation}
    \label{twocases}
a_i<b_k<b_j<a_l\quad\text{ or }\quad b_k<a_i<a_l<b_j.    
  \end{equation}
Consider the first of these two cases. We claim that if it were to
take place, then the configuration 
\[ \tc=\left[(-\infty,b_k)\cup [C\cap(b_k,b_j)]\right]-b_k \]
would be a counterexample to $\kd \delta$.  As $\tc$ has fewer
intervals than $C$, this would contradict the minimality of $C$.

Indeed, consider first a black endpoint $p$ of $C$ between $b_k$ and
$b_j$: $b_k\leq p\leq b_j$, $p\in\B$. We can conclude from the
definition of $J_\B$ that $p+\mu(p)\leq b_j$. Then clearly 
\[ \lambda(\tc|I_{\mu(p)}(p))\geq
\lambda(C|I_{\mu(p)}(p))\geq1-\delta.\]
The proof is analogous when $b_k\leq p\leq b_j$ and $p\in\W$.

The second case of \eqref{twocases} is symmetric to the first one. In
this case 
\[ \tc = a_l-[(a_l,\infty)\cup[C\cap(a_i,a_l)]],\]
and the argument is the same as above.
\end{proof}
\begin{cor}
  Either $I_\circ\subset\fb$ or $I_\circ\subset\fw$.
\end{cor}
This immediately follows from Lemmas \ref{icircin} and
\ref{mainlemma}: if an interval $I$ is contained in the union of a
system of intervals, whose any two elements are either disjoint or one
contains the other, then, in fact, $I$ is already contained in one of
the intervals of the system.

Now we are ready to finish the proof of Proposition \ref{fourth}. 
Because of the symmetry of the problem, without loss of generality, we
can assume that $I_\circ\subset J$, where the interval $J$ is a
connected component of $\fb$. By our construction, the interval $J$ is
a subset of $(0,1)$, and it is a union of intervals of the form
$(a_i,b_i)$ and $I_{\mu(p)}(p)$ with $p\in\B$. Thus it satisfies the
conditions of Lemma \ref{first}, and we can conclude that
\[ \lambda(C\cap J) \geq \frac{1-\delta}{1+\delta}|J|.\] 
As $|I_\circ|\leq|J|$ and $\lambda(C\cap J)\leq\lambda(C\cap(0,1))$,
this implies the statement of the Proposition, and the proof is
complete.
\end{proof}

To prove our main Theorem, all that is left is to make a little
calculation.

According to Lemma \ref{third}, we have 
\begin{equation}
  \label{onerho}
  1-\rho\leq 2\delta(1-v_\B)\quad\text{and} \quad\rho\leq2\delta v_\W.
\end{equation}
Adding up the two inequalities we obtain $1\leq2\delta(1+v_\W-v_\B)$, which
can also be written as
\begin{equation}
  \label{penult}
|I_\circ| \geq \frac1{2\delta}-1.  
\end{equation}
In addition, the second inequality of \eqref{onerho} implies that
\begin{equation}
  \label{small}
\rho\leq2\delta.  
\end{equation}

Substituting \eqref{penult} and \eqref{small} into the inequality of
Proposition \ref{fourth}, we obtain
\[ 2\delta\geq\frac{1-\delta}{1+\delta}\left(\frac1{2\delta}-1\right).\]
Expanding this inequality leads to 
\[ 4\delta^3+2\delta^2+3\delta\geq1,\]
which completes the proof of the Theorem.

\begin{remark}\label{last}
  If we could replace $\frac{1-\delta}{1+\delta}$ by
  $\frac1{1+2\delta}$ in Proposition \ref{third}, then the same
  calculation would lead to the inequality
  $8\delta^3+4\delta^2+2\delta\geq1$. This  does not seem
  impossible, because in Lemma \ref{first} we did not use that we are
  dealing with a special system of intervals. This would confirm our
  conjecture, made at the end of \S\ref{upper}.
\end{remark}

\end{document}